\documentclass[11pt]{article}
\usepackage{palatino}
\usepackage{ifthen}
\usepackage{graphicx}
\usepackage{latexsym}
\usepackage{amssymb}
\usepackage{verbatim}
\newcommand{\diam}{{\sf diam}}
\newcommand{\dist}{{\sf dist}}
\newcommand{\ofpn}{{\sf ofc}}
\newcommand{\dg}{{\sf deg}}

\newcommand{\Bx}{\hfill\ensuremath{\Box}}

\def\br{{\mathbf r}}
\def\mt{{\emptyset}}

\providecommand{\binom}[2]{{#1\choose#2}}

\renewcommand{\baselinestretch}{2.0}

\newtheorem{theorem}{Theorem}[section]
\newtheorem{lemma}[theorem]{Lemma}
\newtheorem{proposition}[theorem]{Proposition}
\newtheorem{corollary}[theorem]{Corollary}
\newtheorem{conjecture}[theorem]{Conjecture}

\newtheorem{fact}[theorem]{Fact}

\begin{document}

\title{$t$-Pebbling and Extensions}
\author{David S. Herscovici\thanks{
Department of Mathematics and Computer Science,
CL-AC3,
Quinnipiac University,
275 Mount Carmel Avenue,
Hamden, CT 06518,
\texttt{David.Herscovici@quinnipiac.edu}}\\
Benjamin D. Hester\thanks{
Department of Mathematics and Statistics,
Arizona State University,
Tempe, AZ, 85287,
\texttt{Benjamin.Hester@asu.edu}} \\
Glenn H. Hurlbert\thanks{
Department of Mathematics and Statistics,
Arizona State University,
Tempe, AZ, 85287,
\texttt{hurlbert@asu.edu}}
}
\renewcommand{\baselinestretch}{1.0}
\maketitle
\renewcommand{\baselinestretch}{2.0}

\begin{abstract}

Graph pebbling is the study of moving discrete pebbles from certain
initial distributions on the vertices of a graph to various target
distributions via pebbling moves.  A pebbling move removes two pebbles
from a vertex and places one pebble on one of its neighbors (losing the
other as a toll).  For $t\ge 1$ the $t$-pebbling number of a graph is the
minimum number of pebbles necessary so that from any initial distribution
of them it is possible to move $t$ pebbles to any vertex.

We provide the best possible upper bound on the $t$-pebbling number
of a diameter two graph, proving a conjecture of Curtis, et al.,
in the process.  We also give a linear time (in the number of edges)
algorithm to $t$-pebble such graphs, as well as a quartic time (in the
number of vertices) algorithm to compute the pebbling number of such
graphs, improving the best known result of Bekmetjev and Cusack.

Furthermore, we show that, for complete graphs, cycles, trees, and cubes,
we can allow the target to be any distribution of $t$ pebbles without
increasing the corresponding $t$-pebbling numbers; we conjecture that
this behavior holds for all graphs.

Finally, we explore fractional and optimal fractional versions of
pebbling, proving the fractional pebbling number conjecture of Hurlbert
and using linear optimization to reveal results on the optimal fractional
pebbling number of vertex-transitive graphs.

\end{abstract}

\section{Introduction}

For a graph $G=(V,E)$, a function $D : V \rightarrow \mathbb{N}$ is
called a \emph{distribution on the vertices of $G$}, or a
\emph{distribution on $G$}.  We usually imagine that $D(v)$ pebbles
are placed on $v$ for each vertex $v \in V$.  Let $|D|$ denote the
\emph{size} of $D$, i.e.\ $|D| = \displaystyle{\sum_{v \in V} D(v)}$.
For two distributions $D$ and $D'$ on $G$, we say that $D$
\emph{contains} $D'$ if $D'(v)\leq D(v)$ for all $v\in V$.  We write
$\diam(G)$ for the diameter of $G$ and $\dist(v, w)$ for the distance
from $v$ to $w$ in $G$.  We use $u\sim v$ to denote that $(u,v)\in E(G)$
($u$ and $v$ are \emph{neighbors}) and define $\dg_X(v)$ to be the 
number of neighbors of $v$ in the set $X$.
In addition, we write $v\sim X$ when $\dg_X(v)\ge 1$.
Here $n$ will represent the number of vertices of $G$.  
The following definition stipulates how pebbles
can be transferred from one vertex to another. \\
\textbf{Definition}: A \emph{pebbling move} in $G$ takes two pebbles
from a vertex $v \in V$, which contains at least two pebbles, and
places a pebble on a neighbor of $v$.  Thus, one pebble is lost. \\
For two distributions $D$ and $D'$, we say that $D'$ is
\emph{reachable} from $D$ if there is some (possibly empty) sequence
of pebbling moves beginning with $D$ and resulting in a distribution
which contains $D'$.  We say that the \emph{cost} of such a sequence
is the sum of the number of pebbles lost along the way and the number
of pebbles placed onto $D'$. \\
\textbf{Definitions}: For an integer $t\geq 1$, we say a distribution
$D$ on a graph $G$ is \emph{$t$-fold solvable} if every distribution
with $t$ pebbles on a single vertex is reachable from $D$.  If $t=1$
we say the distribution is \emph{solvable}; otherwise it is
\emph{unsolvable}.  The \emph{$t$-pebbling number} of a graph $G$,
denoted $\pi_t(G)$, is the smallest integer $k$ such that every
distribution $D$ with $|D| \geq k$ is $t$-fold solvable.  The
\emph{pebbling number} of $G$ is $\pi_1(G)$, and we denote it
$\pi(G)$. \\
In a sequence of pebbling moves, a distribution we are attempting to
reach is called a \emph{target} and a vertex we are attempting to
reach is called a \emph{target vertex} or a \emph{root}.  
For a root $\br$, the quantity $\pi_t(G,\br)$ is the smallest integer $k$ 
such that the distribution with $t$ pebbles on $\br$ and $0$ pebbles 
on every other vertex is reachable from every distribution $D$ with 
$|D| \geq k$.

Suppose that, instead of considering all possible distributions of a
given size, we desire the smallest $t$-fold solvable distribution.  In
this spirit, we give the definition of the optimal $t$-pebbling number
of a graph. \\
\textbf{Definition}: For an integer $t\geq 1$, the \emph{optimal
$t$-pebbling number} of a graph $G$, denoted $\pi_t^*(G)$, is the
smallest integer $k$ such that there exists a $t$-fold solvable
distribution $D$ of pebbles on $V$ with $|D|=k$.  The \emph{optimal
pebbling number} of $G$ is $\pi_1^*(G)$, and we denote it $\pi^*(G)$.

We now outline the remainder of the paper.  The main result is
Theorem~\ref{diam2thm}, which provides the best possible upper bound 
$\pi_t(G)\le\pi(G)+4t-4$ when $G$ has diameter $2$, and proves a 
conjecture of \cite{CHHM} as a corollary.
Furthermore, we obtain an algorithm that places $t$ pebbles on any root
from a distribution of $\pi(G)+4t-4$ pebbles on the $n$ vertices of
such $G$, and that runs in at most $6n+\min\{3t,m\}$ steps, where
$G$ has $m$ edges.
We use this to build another algorithm that calculates $\pi(G)$ 
(distinguishing the two cases $n$ and $n+1$) of such $G$ in $O(n^4)$ time,
besting the work of \cite{Bekmet} when $m\gg n$.
Motivated by prior work of Bukh~\cite{Bukh} and Postle
\emph{et~al.}~\cite{diameter 3 and 4}, we consider graphs of larger
diameter at the end of Section~\ref{t pebbling}.  In particular, we
address a conjecture from~\cite{diameter 3 and 4} and show that any
upper bound on the maximum pebbling numbers of such graphs must be at
least exponential in the diameter.

In Section~\ref{extensions}, we consider extensions and
generalizations of $t$-pebbling numbers.  
In the definition of $\pi_t(G)$ the target is any distribution of $t$
pebbles that all sit on the same vertex.
In the definition of $\pi(G,t)$ the target is any distribution of $t$
pebbles whatsoever.
Necessarily, $\pi(G,t)\ge\pi_t(G)$.
We prove that $\pi(G,t)=\pi_t(G)$ when $G$ is a complete graph, cycle, tree,
or cube, and we conjecture that equality holds for all $G$.
In addition, we prove
a conjecture of \cite{hurlbert URL} (Theorem~\ref{fractional pebbling}), 
which states that the fractional
pebbling number of a graph $G$ is $2^{\diam (G)}$.  In exploring these
extensions, we are naturally led to consider optimal fractional
pebbling numbers, which provide some combinatorial insight into the
fractional world of pebbling.  In addressing fractional optimal
pebbling numbers, we see that they can be found by appealing to linear
optimization.  In the end of the paper, we exploit this fact to
discover results concerning optimal fractional pebbling numbers of
certain graphs and classes of graphs.

\section{$t$-Pebbling}\label{t pebbling}

In Section~\ref{TC}, we describe prior work on pebbling in trees and cycles.  
In Section~\ref{DBB}, we prove a bound on $\pi_t(G,\br)$ which will be useful 
in subsequent sections.  In Section~\ref{GD2}, we prove our main result, 
which provides an upper bound on the $t$-pebbling number of a graph which 
has diameter $2$.  In Section~\ref{GLD}, we address graphs of larger diameter 
by strengthening the known asymptotic lower bound on the value $\pi(n,d)$, 
the maximum pebbling number of an $n$-vertex graph with diameter $d$.

\subsection{Trees and Cycles}\label{TC}

To find the pebbling number of a vertex $\br$ in a tree $T$,
Chung~\cite{Hypercubes} defined $T_\br^*$ as the directed graph in which
all edges in $T$ are directed toward $\br$.  She then described
\emph{path partitions} and \emph{maximal path partitions} in the tree
$T_r^*$, which we generalize to describe path partitions in the
undirected tree $T$ as well. \\
\textbf{Definition} (Chung~\cite{Hypercubes}): A \emph{path partition} of a
undirected tree $T$ or of a tree $T_\br^*$ in which all edges are
directed toward the vertex $\br$ is a partition of the edges of the tree
into sets in such a way that the edges in each set in the partition
form a path in $T$, or a path directed toward $\br$ in $T_\br^*$.  The
\emph{path-size sequence} of a path partition is the sequence of
lengths of the paths in nonincreasing order, $a_1 \geq a_2 \geq \cdots
\geq a_k$.  A \emph{maximal path partition} in $T$ or in $T_\br^*$ is a
path partition whose path-size sequence is lexicographically
greatest. \\
Chung found $\pi_t(T, \br)$, and Bunde \emph{et~al.}~\cite{BCCMW} gave
$\pi_t(T)$.  We present these results as Theorem~\ref{t-pebbling
trees}.
\begin{theorem}[Chung~\cite{Hypercubes}; Bunde \emph{et~al.}~\cite{BCCMW}]
\label{t-pebbling trees}
If $\br$ is a vertex in a tree $T$, then $\pi_t(T,\br)$ is given by
\[ \pi_t(T,\br) = 2^{a_1}t + 2^{a_2} + \ldots + 2^{a_k} - k + 1, \]
where $a_1, a_2, \ldots, a_k$ is the path-size sequence of a maximal
path partition of $T_\br^*$.  Then $\pi_t(T)=\pi_t(T,\br)$, where $\br$ is chosen 
to be the root corresponding to a maximal path partition of $T$.
\end{theorem}
Although it was certainly clear from Chung's work, it appears that no
one has formally stated and proved that moving a pebble to $\br$ costs
at most $2^{a_1}$ pebbles from the rest of the graph.  
We prove this now.
\begin{proposition}
\label{cost of 2^a1 pebbles in tree}
Let $\br$ be any vertex in the tree $T$ and suppose $D$ is a
distribution on $T$ from which $t$ pebbles can be moved to $\br$.  Then
it is possible to move $t$ pebbles to $\br$ at a cost of at most
$2^{a_1} t$ pebbles from the rest of the graph, where $a_1 =
\displaystyle{\max_{v \in V(T)} \dist (\br, v)}$.  In particular, $t$
pebbles can be moved to any vertex at a cost of at most
$2^{\diam(G)}t$ pebbles from the rest of the graph.
\end{proposition}
\textbf{Proof}: Let $S$ be a minimal sequence of pebbling moves that
places $t$ pebbles on $\br$.  For every $i \in \{0, 1, \ldots, a_1 \}$
let $L_i = \{u \in V(G) \mid \dist(\br, u) = i \}$, so $L_i$ is the
$i$th level in the tree rooted at $\br$.  For $i < a_1$ let $n_i$ denote
the number of pebbling moves in $S$ from $L_{i+1}$ to $L_i$.  Now $n_0
= t$, and for larger $i$ we need at most $2n_{i-1}$ moves onto $L_i$
to make $n_{i-1}$ moves onto $L_{i-1}$; therefore, by induction, we
have $n_i \leq 2^i t$.  Thus, the number of pebbling moves in $S$ is
at most $\displaystyle{\sum_{i=0}^{a_1-1} 2^i t = (2^{a_1}-1) t =
2^{a_1} t - t}$.  Each such pebbling move results in the loss of a
pebble.  Thus, along with the $t$ pebbles that ends up on $\br$, at most
$2^{a_1} t$ pebbles are removed from the rest of the graph.~\Bx

The $t$-pebbling number $\pi_t(C_n)$ is given in~\cite{cycles}.  
We will refer to this result in Section~\ref{extensions}.
%
%
\begin{proposition}[Herscovici~\cite{cycles}]
\label{t pebbling cycles}
The $t$-pebbling number of a cycle is given by
\[ \begin{array}{c}
\pi_t(C_{2k}) = 2^k \cdot t \\
\pi_t(C_{2k+1}) = \frac{2^{k+2} - (-1)^k}{3} + 2^k(t-1).
\end{array} \]
\end{proposition}

\subsection{A Distance-Based Bound}\label{DBB}
In this section, we prove Theorem~\ref{radius bound}.
\begin{theorem}
\label{radius bound}
Suppose $\br$ is a vertex in a graph $G$ with the property that
$\dist(v,\br) \leq d$ for every vertex $v$ in $G$.  Then
\begin{equation}
 \pi_t(G,\br) \leq \frac{2^d-1}{d}(n-1) + 2^d (t-1) + 1.
\label{bound}
\end{equation}
Furthermore, if there are at least $\frac{2^d-1}{d}(n-1) + 1$ pebbles
on the graph, moving a pebble to $\br$ costs at most $2^d$ pebbles from
the rest of the graph.
\end{theorem}
%
%
%
%
%
\textbf{Proof of Theorem~\ref{radius bound}}: Let $T$ be a spanning
tree of $G$ obtained by doing a breadth-first search from $\br$.  Since
$T$ is a spanning subgraph of $G$, we have $\pi_t(G,\br) \leq
\pi_t(T,\br)$.  Because $T$ was obtained by a breadth-first search from
$\br$, we have $\dist_G(\br,v) = \dist_T(\br,v)$ for every vertex $v$ in
$G$.  Therefore, it suffices to show that~(\ref{bound}) holds when $G
= T$.  Let $a_1, a_2, \ldots, a_k$ be the path-size sequence of a
maximal path partition of $T_\br^*$.  In particular, we have each $a_i
\leq d$.  Then we have
\begin{equation}
\label{sum of a_i}
\sum_{i=1}^k a_i = |E(T)| = |V(T)| - 1 = n-1,
\end{equation}
since each edge in the tree is in exactly one part in the partition.
From Theorem~\ref{t-pebbling trees}, we also have
\[ \pi_t(T,\br) = 2^{a_1}t + 2^{a_2} + \ldots + 2^{a_k} - k + 1, \]
which, we can rewrite as
\[
\pi_t(T,r) = \sum_{i=1}^k (2^{a_i} - 1) + 2^{a_1}(t-1) + 1 =
\sum_{i=1}^k \left[ \left( \frac{2^{a_i} - 1}{a_i}\right) a_i \right]
+ 2^{a_1}(t-1) + 1.
\]
Now since each $a_i$ is at most $d$, we use the fact that
$(2^d-1)/d$ is an increasing function on the positive integers
to obtain
\[
\pi_t(T,\br) \leq \sum_{i=1}^k \left[ \left( \frac{2^d - 1}{d} \right)
  a_i \right] + 2^d(t-1) + 1 = \left( \frac{2^d - 1}{d} \right)
\sum_{i=1}^k a_i + 2^d(t-1) + 1.
\]
But from~(\ref{sum of a_i}), we find
\[ \pi_t(G,\br) \leq \pi_t(T,\br) \leq \frac{2^d - 1}{d} (n-1) + 2^d(t-1)
+ 1. \]
Finally, since each move is made along a directed edge in $T_\br^*$, by
Proposition~\ref{cost of 2^a1 pebbles in tree}, at most $2^d$ pebbles
from the rest of the graph are consumed.~\Bx

Curtis \emph{et~al.}\ proved Theorem~\ref{n+7t-6 diameter 2}:
\begin{theorem}[Curtis \emph{et\ al.}~\cite{CHHM}]
\label{n+7t-6 diameter 2}
For any integer $t \geq 1$, if $G$ is a graph with diameter $2$, then
$\pi_t(G)\leq n+7t-6$.
\end{theorem}
The proof of Theorem~\ref{n+7t-6 diameter 2} can be generalized to
prove Theorem~\ref{n+7t-6 radius 2}.
\begin{theorem}[Curtis \emph{et\ al.}~\cite{CHHM}]
\label{n+7t-6 radius 2}
If $\br$ is a vertex in $G$ such that $\dist(\br, v) \leq 2$ for every
vertex $v$ in $G$, then $\pi_t(G, \br)\leq n+7t-6$.
\end{theorem}
Using Theorem~\ref{radius bound} with $d=2$ gives the bound $\pi_t(G,
\br) \leq 1.5n + 4t - 4.5$.  Thus, Theorem~\ref{radius bound} represents
an improved bound on that given by Theorem~\ref{n+7t-6 radius 2} when
$t > \frac{n+3}{6}$.  In Section~\ref{diameter 2 section}, we further
improve this bound when the diameter of the graph is $2$.

\subsection{Graphs of Diameter~2}\label{GD2}
\label{diameter 2 section}

We prove Theorem~\ref{diam2thm}, which gives a bound on the
$t$-pebbling number of graphs with diameter $2$.
\begin{theorem}
\label{diam2thm}
If $G$ is a graph with diameter $2$ then $\pi_t(G) \leq \pi(G)+4t-4$.
\end{theorem}

Star graphs, denoted $K_{1,p}$, feature prominently in our proof, so we
define them now. \\
\textbf{Definition}: If $p \geq 2$, the \emph{star on $p+1$ vertices},
denoted $K_{1,p}$, is the graph whose vertex and edge sets are given by
$V(K_{1,p}) = \{ u, v_1, v_2, \ldots, v_p \}$ and $E(K_{1,p}) = \{ (u,
v_i) : 1 \leq i \leq p \}$.  We call $u$ the \emph{center} of the star,
and we call the $v_i$'s its \emph{leaves}.  By abuse of notation,
we identify the vertex set $V$ with the star $K_{1,p}$ if $|V| = p+1$
and the subgraph induced by $V$ contains $K_{1,p}$.

To prove Theorem~\ref{diam2thm}, we let $D$ be a distribution of
$\pi(G)+4t-4$ pebbles on $G$ for some $t \geq 2$ (there is nothing to
show if $t=1$), and we show that $t$ pebbles can be moved to the
vertex $\br$.  We assume by induction that $\pi_{t-1}(G) \leq \pi(G) +
4(t-1)-4 = \pi(G) + 4t-8$.  Therefore, if we could move a pebble to
$\br$ at a cost of no more than four pebbles, we could use the remaining
$\pi(G) + 4t-8$ pebbles to put $t-1$ additional pebbles on $\br$.  We
show that if putting a pebble on $\br$ requires using five
pebbles, then $n+4t-4$ pebbles are sufficient to put $t$ pebbles on
$\br$.  To do this, we note that if four pebbles are not sufficient to
move a pebble onto $\br$, this places certain constraints on $D$.
Lemmas~\ref{distance 2 from r}, \ref{middle unoccupied},
and~\ref{middle distinct} formalize this idea.
\begin{lemma}
\label{distance 2 from r}
Suppose $G$ is a graph with diameter $2$, and $D$ is a distribution
on $G$ from which any sequence of pebbling moves that puts a pebble on
the vertex $\br$ requires at least five pebbles.  Then every vertex has
at most three pebbles, and no vertex with two or three pebbles can be
adjacent to $\br$.
\end{lemma}
\begin{lemma}
\label{middle unoccupied}
Suppose $G$ is a graph with diameter $2$, and $D$ is a distribution
that satisfies the condition of Lemma~\ref{distance 2 from r}.  Let
$v_i$ be a vertex with at least two pebbles.  Then there is a vertex
$w_i$ adjacent to both $v_i$ and $\br$.  Furthermore, every such $w_i$
is unoccupied.
\end{lemma}
\begin{lemma}
\label{middle distinct}
Suppose $G$ is a graph with diameter $2$, and $D$ is a distribution
that satisfies the condition of Lemma~\ref{distance 2 from r}.
Suppose further that $v_i$ and $v_j$ are distinct vertices which each
have two or three pebbles.  Then the vertices $w_i$ and $w_j$ from
Lemma~\ref{middle unoccupied} are also distinct.
\end{lemma}
\textbf{Proof of Lemmas~\ref{distance 2 from r}, \ref{middle
    unoccupied}, and~\ref{middle distinct}}: Let $v_i$ be any vertex
with at least two pebbles.  We note that if $v_i$ were adjacent to $\br$
(or if $v_i = \br$), two pebbles would be sufficient to reach $\br$.
Since $\diam(G) = 2$, there is a vertex $w_i$ that is adjacent to both
$v_i$ and to $\br$.  Therefore, four pebbles on $v_i$ would be
sufficient to reach $\br$ (completing the proof of Lemma~\ref{distance 2
  from r}), and if $w_i$ were occupied, two pebbles on $v_i$ and one
pebble on $w_i$ would be sufficient to reach $\br$ (establishing
Lemma~\ref{middle unoccupied}).  Finally, if any $w_i$ were adjacent
to both $v_i$ and $v_j$ in Lemma~\ref{middle distinct}, then we could
move one pebble onto $w_i$ from $v_i$ and another from $v_j$, and from
there we could move a pebble onto $\br$ at a total cost of four
pebbles.~$\Bx$


Before proving Theorem \ref{diam2thm}, we introduce the new concept of
resolving a subgraph.  Essentially, if $D$ is a distribution on a graph
$G$ and $H$ is an edge subgraph of $G$ (i.e. $E(H)\subseteq E(G)$), then
we define the distribution $D_H$ on $G$ by $D_H(v)=D(v)$ for $v\in H$ and
0 otherwise.  Vaguely, for some root $\br$ of $G$, when we say to resolve
$H$, we mean to place as many pebbles on $\br$ as possible from the
distribution $D_H$.  When the time comes, for certain subgraphs having
particular distributions, we will remove ambiguity by describing the
necessary pebbling steps in sufficient detail.

For example, let $H$ be the star $K_{1,p}$ with $p\ge 2$, having 3 pebbles
on each of its leaves and at most 2 pebbles on its center.  In this case,
we \emph{resolve the star} by, first, moving pebbles from some of its
leaves through the center and onto $l$ other leaves $L$ so that every leaf
has at most 4 pebbles and so that $l$ is maximized and, second, moving
$l$ pebbles onto $\br$ from $L$ (which is possible because $\diam(G)=2$).

\begin{lemma}
\label{consuming n3s}
Let $H$ be a star $K_{1,p}$ with $p \geq 2$, whose center vertex has
$i \leq 2$ pebbles and whose leaves each have three pebbles.
Then resolving $H$ puts $l=\lfloor (p+i)/3 \rfloor$
pebbles on $\br$.  Moreover, for $l'=(p+i) \bmod 3$ (so that $p+i = 3l+l'$
and $0 \leq l' \leq 2$), there remain $l'$ leaves with 3 pebbles each,
and the number of leaves used in the resolution equals $3l-i$.
\end{lemma}
\textbf{Proof}: 
For every three leaves we can move pebbles from two leaves to the center
and then one pebble to the third.  For every pebble already on the center
we save a pebbling step from a leaf to the center.
This uses $3l-i$ outer vertices, so the number of unused outer vertices 
is $p-(3l-i) = p+i-3l = l'$.~$\Bx$

Notice that Lemma \ref{consuming n3s} does not hold for $p=1$.
For this reason, our strategy for placing $t$ pebbles on $\br$ in $G$ 
will also use matchings between vertices having 3 pebbles.
We \emph{resolve a matching} $M$ such as this by, first, arbitrarily moving 
one pebble across each edge and, second, moving a pebble from each
recipient vertex to $\br$.  The number of pebbles that reach $\br$ equals
the number $m$ of edges of $M$, and all $2m$ vertices of $M$ are used in
this resolution.

We are now ready to prove Theorem~\ref{diam2thm}. \\
\textbf{Proof of Theorem~\ref{diam2thm}}: 
As discussed above, we may assume the hypothesis of Lemma 
\ref{distance 2 from r}; otherwise induction suffices.
Our proof is algorithmic.
In the first stage we resolve a matching; in the second we iteratively
resolve stars.
By some careful counting arguments we will show that $t$ pebbles reach $\br$.

We begin by defining, for $0\le k\le 3$, the set $V_k=\{v\mid D(v)=k\}$, 
with $n_k=|V_k|$.
Let $M$ be a maximal matching in the subgraph $G[V_3]$ induced by the
vertices of $V_3$, and denote its number of edges by $m$.
Resolve $M$ and let $D'$ be the resulting distribution.
For $0\le k\le 3$ define $V_k'=\{v\in V-V(M)\mid D'(v)=k\}$.
Notice that $V_3'$ is independent in $G$.

Next we initialize the sets $S_k=\mt$ for $0\le k\le 2$ and 
$L=\mt$, and iterate the following steps.
For a vertex $v$ define $d_3'(v)=\dg_{V_3'}(v)$.
Now for any $k$ choose some $v\in V_k'$ with $d_3'(v)\ge 3-\lfloor k/2\rfloor$,
if one exists (necessarily $k\le 2$), and let $S$ be the star consisting of 
the center $v$ and all its neighbors in $V_3'$.
Resolve $S$, put $v$ in $S_k$, add the leaves of $S$ to $L$,
redefine the notation $D'$ for the 
resulting distribution, and likewise redefine the sets $V_k'$ accordingly.
At some point the algorithm halts because no possibilities remain for
choosing an appropriate center $v$.
Write $s$ for the number of pebbles that the stars contribute to $\br$.

Over the course of the algorithm, Lemma \ref{consuming n3s} implies that
\begin{equation}\label{pebbles solved}
m+s=m+\sum_{k=0}^2\sum_{v\in S_k}l(v)
\end{equation}
pebbles have reached $\br$, where $l(v)=\lfloor(d_3'(v)+k)/3\rfloor$, 
and that 
\begin{equation}\label{vertices used}
2m+\sum_{k=0}^2\sum_{v\in S_k}(3l(v)-k)=2m+3s-s_1-2s_2
\end{equation}
vertices from $V_3$ were used to send pebbles in the process, where $s_k = |S_k|$.
The proof will be complete when we show that $m+s\ge t$.

At this time set $U=V_3'-L$, $W=\{v\mid\dg_U(v)\ge 2\}$,
and $Y=\{v\mid \br\sim v\sim V_3\}$.
We note that $u=|U|=n_3-2m-3s+s_1+2s_2$ (by (\ref{vertices used})) and
$|W|\ge \binom{u}{2}\ge u-1$ (since $\diam(G)=2$ and no more stars exist).
Observe also that $V_3\cap Y=W\cap Y=\mt$ because of the cost 5 assumption.
Hence, for $Z=V_3\cup W\cup Y$ we compute
\begin{eqnarray*}
z=|Z|&=&|V_3|+|Y|+|W-V_3|\\
&\ge&n_3+n_3+(u-1-m)\\
&=&3n_3-3m-3s+s_1+2s_2-1,
\end{eqnarray*}
the appearance of $m$ coming from the fact that, if $v\in W\cap V_3$,
then having two neighbors in $U$ when the algorithm halts means that
$D'(v)=0$, $v\notin\cup_{k=0}^2S_k$, and consequently that $v$ was a
recipient in the resolution of $M$.  The number of such vertices was
exactly $m$.  We also can calculate the number of pebbles originally
on $Z$.
\begin{eqnarray*}
|D(Z)|&=&|D(V_3)|+|D(Y)|+|D(W-V_3)|\\
&\le&3n_3+0+s_1+2s_2\\
&\le&z+3m+3s+1.
\end{eqnarray*}

To complete the analysis we set $X_k=V_k-Z$ ($0\le k\le 2$) and
$x_k=|X_k|$; then $n=z+x_0+x_1+x_2$.  Now we have
\begin{eqnarray*}
z+x_0+x_1+x_2+4t-4&=&n+4t-4\\
&\le&|D|\\
&=&|D(Z)|+x_1+2x_2\\
&\le&z+3m+3s+1+x_1+2x_2.
\end{eqnarray*}
Combined with the fact that $x_0\ge x_2+1$ (since 
$X_0 \supseteq \{\br\}\cup\{v\mid\br\sim v\sim X_2\}$), this implies that
$4t-4\le 3m+3s$, from which follows
\begin{eqnarray*}
m+s&\ge&\left\lceil\frac{4t-4}{3}\right\rceil\\
&=&(t-1)+\left\lceil\frac{t-1}{3}\right\rceil\\
&\ge&t,
\end{eqnarray*}
since $t\ge 2$.
This completes the proof.~$\Bx$


Pachter, Snevily, and Voxman~\cite{PSV} found the pebbling number for
graphs with diameter~2.  In particular, they showed the following theorem.
\begin{theorem}[Pachter \emph{et\ al.}~\cite{PSV}]
\label{diam2bound}
If $G$ is a graph with diameter $2$, then $\pi(G)\leq n+1$.
\end{theorem}

Putting Theorems~\ref{diam2thm} and~\ref{diam2bound} together, gives
us Corollary~\ref{n+4t-3}, first conjectured in~\cite{CHHM}.
\begin{corollary}
\label{n+4t-3}
For any integer $t \geq 1$, if $G$ is a graph with diameter
$2$, then $\pi_t(G) \leq n+4t-3$.~$\Bx$
\end{corollary}

One might ask whether there are any diameter~$2$ graphs $G$ and any
values of $t$ where the inequality in Theorem~\ref{diam2thm} is
strict, i.e.\ $\pi_t(G) < \pi(G) + 4t-4$.  Indeed there are.
Proposition~\ref{Kn-e} shows that the difference can be as large as
$\pi(G)-4$.
\begin{proposition}
\label{Kn-e}
For $n \geq 3$, let $G_n$ be the graph obtained by removing the edge
$\{ v_1, v_2 \}$ from the complete graph $K_n$.  Then for any $t \geq
n-2$, we have $\pi_t(G_n) = 4t$.
\end{proposition}
\textbf{Proof}: We have $\pi_t(G) \geq 4t$, since placing $4t-1$
pebbles on $v_2$ creates a distribution from which $t$ pebbles cannot
be moved to $v_1$, so we need to show $\pi_t(G) \leq 4t$.  We use
induction on $n$, and later, induction on $t$ as well.  The basis is
$n=3$.  Now $G_3$ is the path on the vertices $\{ v_1, v_3, v_2 \}$,
in that order, so $\pi_t(G_3) = 4t$, as desired.

For $n > 3$, we assume $\pi_{t'}(G_{n-1}) = 4t'$ whenever $t' \geq
n-3$.  We also assume $v_1$ is the target (or $v_2$); otherwise,
applying Theorem~\ref{radius bound} with $d=1$ gives $\pi_t(G) \leq n
+ 2t - 2 \leq 3t$.

Let $D$ be a distribution of $4t$ pebbles on $G_n$, and let $p_i =
D(v_i)$.  We assume without loss of generality that $p_3 \geq p_4 \geq
\cdots \geq p_n$.  The rest of our argument depends on the value of
$p_n$.  If $p_n = 0$ we have $4t$ pebbles on the vertices $\{ v_1,
v_2, \ldots, v_{n-1} \}$.  These vertices induce a subgraph isomorphic
to $G_{n-1}$, so by our inductive assumption, $4t$ pebbles are
sufficient to $t$-pebble $v_1$.  If $p_n = 1$, we note that some
vertex has two pebbles, and since $v_n$ is adjacent to every other
vertex, we can put a second pebble onto $v_n$, and from there, we can
put a pebble on $v_1$.  After that, the remaining $4t-3$ pebbles on
$\{ v_1, v_2, \ldots, v_{n-1} \}$ suffice to put an additional $(t-1)$
pebbles on $v_1$, again by our inductive assumption.

Finally, if $p_n \geq 2$, we resort to induction on $t$.  We note that
each of the $(n-2)$ vertices in $\{ v_3, v_4, \ldots, v_n \}$ has at
least two pebbles and is adjacent to $v_1$.  Therefore, at least
$(n-2)$ pebbles can be moved to $v_1$.  Thus, if $t = n-2$, we are
done.  Otherwise, $t \geq n-1$, so we move one pebble from $v_n$ to
$v_1$.  Since $t-1 \geq n-2$, we may assume by induction on $t$ that
the remaining $4t-2$ pebbles are sufficient to move $(t-1)$ additional
pebbles onto $v_1$.~$\Bx$

\subsection{Algorithmic Results for Diameter 2 Graphs}\label{new algos}

We remark that, for $|D|\ge \pi(G)+4t-4$, we can place $t$ pebbles on
any root $\br$ in time that is linear in $n$ and $t$: at most $6n+3t$
steps are required.  To see this, if one can place a pebble on $\br$
with cost at most 4, a breadth-first search from $\br$ will reveal it
in at most $n$ steps, since the only possibilities (besides already
having a pebble on $\br$) are (\emph{i}) a path $(\br,u)$ with
distribution $(0,\ge 2)$, (\emph{ii}) a path $(\br,u,v)$ with
distribution $(0,1,\ge 2)$, (\emph{iii}) a path $(\br,u,v,w)$ with
distribution $(0,1,1,\ge 2)$, (\emph{iv}) a path $(\br,u,v)$ with
distribution $(0,0,\ge 4)$, and (\emph{v}) a star $(\br,u,v,w)$ with
center $u$ and distribution $(0,0,\ge 2,\ge 2)$.  In fact, it is not
difficult to see that all these solutions can be found and resolved in
at most $n+3t$ steps, since each resolution involves at most 3 edges
and the breadth-first search does not need to be repeated.

If none of these possibilities exist, then our algorithm will place
the required pebbles on $\br$ in at most $3n$ steps when $t\ge 2$.
Indeed, since the directed graph formed by orienting the edges of $G$
in the direction of pebbling can be seen to be acyclic (which is not
necessarily the case above), the total
number of edges used in the algorithm is at most $n$; consequently
we can implement the algorithm so as to postpone the actual pebbling
steps so that each such edge is traversed once.
Thus we only need to count the number of steps required to find
the matching and stars.
It takes $n$ steps to sort the vertices into the appropriate $V_k$.
Then a maximal matching in $G[V_3]$ can be found in at most $v_3/2$
steps (we are fortunate not to need a maximum matching).
Finally, if we search first for star centers with the most pebbles,
then we don't repeat vertices in our search, so the total time to 
find all stars is at most $n-v_3$, making the number of steps
for finding the matching and all stars at most $n$.
Including sorting and resolution, we use at most $3n$ steps.

As mentioned, in the proof of Theorem \ref{diam2thm}, we only used 
our algorithm when $t\ge 2$.
Here we observe that it actually works with a slight modification
when $t=1$ as well.
The only difference is that we look not only for matching edges
in $V_3$, but also for edges between $V_2$ and $V_3$.
Hence if we have been unsuccessful to this point in placing a
pebble on $\br$, we will show that, to avoid the contradiction 
that $|D|<\pi(G)$, we can find one final solution method in 
at most $2n$ steps.
To do so, define the sets $V_{k,d}=\{v\mid D(v)=k, \dist(v,\br)=d\}$,
for which we know that $0\le d\le 2$ and $0\le k<2^d$ (so that
$V_k=V_{k,2}$ for $k\in\{2,3\}$).
In addition, we know that the neighbors of $V_2\cup V_3$ that lie in 
$V_{0,1}$ are distinct, and that the common neighbors of
pairs of points $(u,v)\in V_3\times (V_2\cup V_3)$ 
lie in $V_{0,2}$ and are distinct also.
Thus, in order that $|D|\ge n$ we must have 
$|V_2|+2|V_3|\ge |V_{0,1}|+|V_{0,2}|+1$;
that is, there must be enough extra pebbles from $V_2\cup V_3$
to compensate for the lack of pebbles on $V_0$, including $\br$.
In particular, this implies that 
$|V_3|>|V_{0,2}|\ge\binom{|V_3|}{2}+|V_3||V_2|$,
which means that $|V_3|\in\{1,2\}$, $|V_2|=0$, 
$|V_{0,2}|=\binom{|V_3|}{2}$, and $\pi(G)\le |D|=n$.
Since we know that it is possible to place a pebble on $\br$,
and we know it cannot begin by moving anything from $V_3$ to
$V_{0,1}$, it must arise from moving pebbles from $S=V_3\cup V_{0,2}$
along a path through $V_1$ to $\br$.
Such a path $P$ can be found from a breadth-first search from $\br$
in $G[V-V_{0,1}]$; whichever vertex $v$ from $S$ is found first
determines whether to move directly along $P$ from $v\in V_3$ 
or first to move 2 pebbles to $v\in V_{0,2}$ and then on to $P$.
Finding and pebbling along $P$ takes at most $2n$ steps.

We mention finally that if $3t>m=m(G)$, the number of edges of $G$, we
can perform the same trick in the first stage above as in the second,
namely, that we postpone the actual resolution until the end, using at
most $m$ steps instead of $3t$.  
We record this in the following theorem.

\begin{theorem}\label{running time}
If $G$ is a graph with $n$ vertices, $m$ edges, and diameter 2, and 
$D$ is a distribution of size $\pi(G)+4t-4$, then $t$ pebbles can be 
placed on any root $\br$ in at most $6n+\min\{3t,m\}$ steps.
\end{theorem}

The only other algorithmic results known for diameter two graphs are
found in \cite{Bekmet}.
There the authors consider the case when $t=1$ and $|D|<\pi(G)$, and
present algorithms that determine the solvability of $D$ in
polynomial time on graphs of constant bounded connectivity, among
other cases.
Their algorithm uses the characterization found in \cite{CHH} (see
also \cite{BS}) for diameter 2 graphs with pebbling number $n+1$, 
as opposed to $n$.
They also use this characterization to give an $O(n^3m)$ algorithm 
for recognizing the difference.  
Here we use the $t=1$ portion of our algorithm to produce an $O(n^4)$ 
algorithm that recognizes the difference, which improves their result
when $m\gg n$.

\begin{theorem}\label{Class char}
If $G$ is a graph with $n$ vertices and diameter 2, then it can be
determined whether $\pi(G)=n$ or $n+1$ in $O(n^4)$ time.
\end{theorem}

\textbf{Proof}: 
We note that $\pi(G)=n+1$ if and only if there is a distribution $D$
of $n$ pebbles that cannot reach some $\br$.
For such a $D$ we have argued that it must have $|V_3|\in\{1,2\}$,
along with the other conditions mentioned above.
We use breadth-first search to determine distances from $\br$ and, 
for each $i\in\{1,2\}$ we pick $i$ vertices at distance 2 from $\br$
each of which has a unique common neighbor with $\br$, determining 
both $V_3$ and $V_{0,1}$ (such available choices can be filtered 
during the breadth-first search).
We ignore the choice unless each chosen vertex has a unique common
neighbor with $\br$, determining $V_{0,1}$, and when $i=2$ the pair
in $V_3$ has a unique common neighbor, determining $V_{0,2}$.
Now $\pi(G)=n$ if and only if there is a path from $V_3\cup V_{0,2}$
to $\br$ in $G[V-V_3-V_0]$.
The total time required to complete these constructions and checks
is at most $n+(n+n^3/2)n=n^4/2+n^2+n$.~$\Bx$

\subsection{Graphs of Larger Diameter}\label{GLD}

It would be interesting to expand our methods with graphs of diameter $2$ to graphs with larger
diameter. Postle, Streib, and Yerger~\cite{diameter 3 and 4} announced
Theorem~\ref{1.5n+2}, strengthening a result of Bukh~\cite{Bukh}.
\begin{theorem}[Postle \emph{et~al.}~\cite{diameter 3 and 4}]
If $G$ is a graph with diameter~$3$, then $\pi(G) \leq 1.5n + 2$.
\label{1.5n+2}
\end{theorem}

Let $\pi(n,d)$ denote the maximum pebbling number of an $n$-vertex
graph which has diameter $d$.  Bukh proved Theorem~\ref{Bukh bound}.
\begin{theorem}[Bukh~\cite{Bukh}]
\label{Bukh bound}
There are constants $c, N, D$ such that, for all $n>N$ and $d>D$, we have
$\pi(n,d) \ge \left( \frac{2^{\left \lceil \frac{d}{2} \right
\rceil} - 1}{\left \lceil \frac{d}{2} \right \rceil} \right) n + c$.
\label{asymptotic LB}
\end{theorem}
Postle \emph{et~al.}\ conjectured the following asymptotic upper bound
on $\pi(n,d)$.

\begin{conjecture}[Postle \emph{et~al.}~\cite{diameter 3 and 4}]
There are constants $C, N$ and a function $f$ on the positive integers
such that, for all $n>N$ we have $\pi(n,d) \le \left( \frac{2^{\left
    \lceil \frac{d}{2} \right \rceil} - 1}{\left \lceil \frac{d}{2}
  \right \rceil} \right) n + Cf(d)$.
\label{general diameter}
\end{conjecture}
We show that if Conjecture~\ref{general diameter} holds, then $f(d)$
is at least exponential in $d$, for all $n$.  We do this by creating
$n$-vertex graphs of diameter $d$ that have large pebbling numbers for
all $n$ and $d$.  Given positive integers $n$ and $d$, we build the
graph $G_{n,d}$ as follows. \\
If $d = 2k$, choose a vertex $v$ and build $\left\lfloor \frac{n-1}{k}
\right\rfloor$ paths of length $k$ beginning at $v$ which are disjoint
(except of course at $v$).  If the number of vertices at this point is
smaller than $n$, add one more path of length $n-k\left\lfloor
\frac{n-1}{k} \right\rfloor -1$ which begins at $v$ and is disjoint
from the rest of the graph.  This is the graph $G_{n,2k}$.\\
If $d = 2k+1$, build a clique $K_m$, where $m = \left\lfloor
\frac{n}{k+1} \right\rfloor$.  From each vertex in this $K_m$, build a
path of length $k$ which is disjoint from the rest of the graph.  If
the number of vertices at this point is smaller than $n$, choose any
$v \in V(K_m)$ and add one more path of length $n-m(k+1)$ which begins
at $v$ and is disjoint from the rest of the graph.  This is the graph
$G_{n,2k+1}$.\\
It is easy to check that $G_{n,2k}$ and $G_{n,2k+1}$ have $n$ vertices
and diameters $2k$ and $2k+1$, respectively.
 \begin{proposition}
\label{lb for pi(G)}
$\pi(G_{n,d}) \geq \left(
\frac{2^{\left \lceil \frac{d}{2} \right \rceil} - 1}{\left \lceil
\frac{d}{2} \right \rceil} \right) n + \left(2^d -
  3\left(2^{\left\lceil \frac{d}{2} \right\rceil}-1\right) \right)$
  for all $n$ and $d$.  In particular, if Conjecture~\ref{general
    diameter} holds, then $f(d)\ge c2^d$ for some $c$ and all 
large enough $d$.
\end{proposition}
\textbf{Proof}: If $d=2k$ for some $k$, then $G_{n,d}$ is a tree.  In
a maximal path partition, there is one path of length $2k$ and there
are $\left\lfloor \frac{n-1}{k} \right\rfloor-2$ paths of length $k$.
Thus, from Theorem \ref{t-pebbling trees}, we have
\begin{eqnarray*}
\pi(G_{n,d}) & \geq & 2^{2k}+2^{k}\left(\left\lfloor
 \frac{n-1}{k} \right\rfloor -2 \right)-\left(\left\lfloor
 \frac{n-1}{k} \right\rfloor-1\right) +1 \nonumber \\
& = & \left( 2^{k} - 1 \right) \left( \left\lfloor
 \frac{n-1}{k} \right\rfloor \right) + \left( 2^{2k} -
 2^{k+1} + 2 \right) \\
& \geq & \left( 2^{k} - 1 \right) \left(
 \frac{n}{k} - 1 \right) + \left( 2^{2k} - 2^{k+1} +
 2 \right) \\
& = & \left( \frac{2^{\left \lceil \frac{d}{2} \right
 \rceil} - 1}{\left \lceil \frac{d}{2} \right \rceil} \right) n +
 \left(2^d - 3\left(2^{\lceil \frac{d}{2}
 \rceil}-1\right)\right).
\end{eqnarray*}
If $d = 2k+1$ for some $k$, we build an unsolvable distribution $D$ on
$G_{n,d}$.  Let the root $\br$ be any leaf which is the endpoint of a
maximum induced path $P$.  Place $2^{2k+1}-1$ pebbles on the other
endpoint of $P$.  Now, for every leaf vertex (disjoint from $P$) that is
distance $k$ from $K_m$, place $2^{k+1}-1$ pebbles.  There are
$\left\lfloor \frac{n}{k+1} \right\rfloor-2$ such vertices.  It is
easy to verify that $D$ cannot send a pebble to $\br$.  Thus, since
$\pi(G_{n,d})\geq|D|+1$,
\begin{eqnarray*}
\pi(G_{n,d}) & \geq & 2^{2k+1} + \left(2^{k+1} - 1 \right)
 \left( \left\lfloor \frac{n}{k+1}
 \right\rfloor -2 \right) \\
& = & \left( 2^{k+1} - 1 \right) \left(
\left\lfloor \frac{n}{k+1}
\right\rfloor \right)+\left(2^{2k+1} - 2^{k+2} + 2
\right) \\
& \geq & \left(2^{k+1} - 1 \right) \left(
 \frac{n}{k+1} - 1 \right) + \left(
 2^{2k+1} - 2^{k+2} + 2 \right) \\
& = & \left( \frac{2^{\left \lceil \frac{d}{2} \right \rceil} -
1}{\left \lceil \frac{d}{2} \right \rceil} \right) n + \left(2^d -
3\left(2^{\lceil \frac{d}{2} \rceil}-1\right)\right),
\end{eqnarray*}
as desired.\Bx

Conjecture~\ref{general diameter} would follow from
Conjecture~\ref{pi(G) <= pi(Gnd)}.
\begin{conjecture}
\label{pi(G) <= pi(Gnd)}
Let $G$ be any $n$-vertex graph with diameter $d$.  Then $\pi(G) \leq
\pi(G_{n,d})$.
\end{conjecture}

\section{Extensions}\label{extensions}

In Section~\ref{ATDtP}, we consider how many pebbles are required to reach
an arbitrary target distribution with $t$ pebbles.  We conjecture an equality 
which would relate $\pi_t(G)$ to a more general pebbling invariant on a graph.  
The truth of the equality would simplify the process of obtaining general results 
about achieving arbitrary target distributions on graphs.  In Section~\ref{FP}, 
we discuss how the $t$-pebbling number of a graph increases as $t$ increases.  
This naturally leads to the discussion of the fractional analogue of the pebbling 
number of a graph.  We show that this value depends only on $\diam(G)$.  
In Section~\ref{OFP}, we analyze the continuous version of optimal pebbling 
and present the corresponding linear optimization problem.

\subsection{Arbitrary Target Distributions with \boldmath{t} Pebbles}\label{ATDtP}
The following definition generalizes the definition of the $t$-pebbling number of a graph. \\
\textbf{Definition}: We define $\pi(G, t)$ as the smallest number of
pebbles such that any target distribution $D$ with $|D| = t$ is
reachable from every distribution $D'$ with $|D'| \geq \pi(G, t)$. \\
Clearly if we can reach any distribution with $t$ pebbles starting
from $D$, we can reach any distribution with $t$ pebbles on a single
vertex.  Therefore, $\pi_t(G) \leq \pi(G, t)$ for every positive
integer $t$.  Conversely, it seems reasonable to believe that if we
have a distribution of pebbles from which we can put $t$ pebbles on
any single vertex of $G$, then any other distribution of $t$ pebbles
is likewise reachable.  For example, if we can put two pebbles either
on the vertex $x$ or the vertex $y$, then we should be able to put one
pebble each on $x$ and $y$.  This suggests the following conjecture.
\begin{conjecture}
For every graph $G$ and every positive integer $t$, we have $\pi(G, t)
= \pi_t(G)$.
\label{targets with t pebbles}
\end{conjecture}
One might be interested in a more general target distribution as a stepping 
stone to some goal.  
Having the equality from Conjecture~\ref{targets with t pebbles} as a tool 
could greatly simplify the necessary analysis.  We prove this conjecture for 
some common graphs.  We start with two lemmas.  
\begin{lemma}
Suppose $G$ is a graph with the property that, for some $t$, whenever
$\pi_{t+1}(G)$ pebbles are on $G$, one pebble can be moved to any
vertex at a cost of at most $\pi_{t+1}(G)-\pi(G, t)$ pebbles.  Then
$\pi(G, t+1)=\pi_{t+1}(G)$.
\label{spending pebbles}
\end{lemma}
\textbf{Proof}: Let $D$ be a distribution on $G$ with $t+1$ pebbles.
Given a distribution of $\pi_{t+1}(G)$ pebbles on $G$, choose one
occupied vertex $v$ in $D$, and spend $\pi_{t+1}(G)-\pi(G, t)$ pebbles to
move a pebble to $v$.  The remaining $\pi(G, t)$ pebbles can be used to
move $t$ additional pebbles to fill out the rest of $D$.\Bx
\begin{lemma}
Suppose $G$ is a graph with the property that for every $t$, if
$\pi_{t+1}(G)$ pebbles are on $G$, one pebble can be moved to any
vertex at a cost of at most $\pi_{t+1}(G)-\pi_t(G)$ pebbles.  Then
$\pi(G,t)=\pi_t(G)$ for all $t$.
\label{pi_t(G) = pi(G, t)}
\end{lemma}
\textbf{Proof}: We use induction on $t$.  When $t=1$ we have $\pi_1(G)
= \pi(G, 1) = \pi(G)$ since the target distributions are the same in
either case.  For larger $t$, if $\pi_t(G) = \pi(G, t)$, then
$\pi_{t+1}(G)-\pi_t(G) = \pi_{t+1}(G)-\pi(G, t)$, so by
Lemma~\ref{spending pebbles}, $\pi_{t+1}(G) = \pi(G, t+1)$. \\
Theorems~\ref{2^d t} and~\ref{Kn T and Cn} gives some classes of
graphs for which Conjecture \ref{targets with t pebbles} holds.
\begin{theorem}
\label{2^d t}
Let $G$ be any graph such that $\pi(G)=2^{\diam(G)}$.  Then for any $t
\geq 1$, $\pi(G,t)=\pi_t(G) = 2^{\diam(G)} t$.  In particular,
Conjecture \ref{targets with t pebbles} holds for complete graphs, 
even cycles, and hypercubes.
\end{theorem}
\textbf{Proof}: By Lemma~\ref{pi_t >= 2^d t}, $2^{\diam(G)} t \leq
\pi_t(G)$.  Conversely, given $2^{\diam(G)} t$ pebbles, we can split
them into $t$ groups of $2^{\diam(G)}$ pebbles each.  Then each group
can be matched to a different pebble in any target distribution with
$t$ pebbles.  Thus, $2^{\diam(G)} t \leq \pi_t(G) \leq \pi(G, t) \leq
2^{\diam(G)} t$, so $\pi_t(G) = \pi(G, t)$.\Bx
\begin{theorem}
\label{Kn T and Cn}
If $G$ is a tree or a cycle, then $\pi(G, t)= \pi_t(G)$.
\end{theorem}
\textbf{Proof}: If $G$ is a tree, by Theorem~\ref{t-pebbling trees}
and Proposition~\ref{cost of 2^a1 pebbles in tree}, the cost of
putting a pebble on any vertex is at most $2^{\diam(G)} = \pi_{t+1}(G)
- \pi_t(G)$.  If $G$ is an even cycle we can apply Theorem~\ref{2^d
  t}, so suppose $G = C_n$ is an odd cycle with vertices $\{ x_1, x_2,
\ldots, x_n \}$ in order with $n = 2k+1$.  We may assume without loss
of generality that $x_n$ is the target vertex.  By Proposition~\ref{t
  pebbling cycles}, $\pi_t(C_n)$ is given by $\pi_t(C_n) =
\frac{2^{k+2} - (-1)^k}{3} + 2^k (t-1)$.
Thus, $\pi_{t+1}(C_n) - \pi_t(C_n) = 2^k$, and $\pi_t(G) \geq 2^{k+1}$
when $t \geq 2$.  In particular, if we have $\pi_{t+1}(G) \geq
2^{k+1}$ pebbles on $C_n$, either we have $2^k$ pebbles on the
vertices $\{ x_n, x_1, x_2, \ldots, x_k \}$ or we have $2^k$ pebbles
on the vertices $\{ x_n, x_{n-1}, x_{n-2}, \ldots, x_{k+1} \}$.  Since
these vertex sets each induce the subgraph $P_{k+1}$, we can move a
pebble to $x_n$ at a cost of at most $\pi(P_{k+1}) = 2^k$ pebbles.\Bx

\subsection{Fractional Pebbling Numbers}\label{FP}
\label{fracpebb}
One might wonder how the $t$-pebbling number of a graph
grows with $t$.  We note that for complete graphs, trees,
cycles, and indeed for all other graphs $G$ for which $\pi_t(G)$ is
known, we have $\pi_{t+1}(G) \leq \pi_t(G) + 2^{\diam(G)}$ for all
$t$.  We raise this observation to the status of a conjecture, and we
prove it for large enough $t$.  Conjecture~\ref{weak diam conj} is a
weaker version of Conjecture~\ref{diamconj}.
\begin{conjecture}
\label{diamconj}
For every graph $G$ and for every $t \geq 1$, we have $\pi_{t+1}(G)
\leq \pi_t(G) + 2^{\diam(G)}$.
\end{conjecture}
\begin{conjecture}\label{weak diam conj}
For every graph $G$ and for every $t \geq 1$, we have $\pi_t(G)
\leq \pi(G) + 2^{\diam(G)} (t-1)$.
\end{conjecture}
Theorem~\ref{diam2thm} proves Conjecture~\ref{weak diam conj} for all
graphs with diameter~$2$.  Combining Conjecture~\ref{weak diam conj}
with Theorem~\ref{1.5n+2} gives us Conjecture~\ref{diam 3 t pebbling
conjecture}, and combining it with Conjecture~\ref{general diameter}
gives Conjecture~\ref{t pebbling general diam}.
\begin{conjecture}
\label{diam 3 t pebbling conjecture}
If $G$ is a graph with diameter~$3$, then $\pi_t(G) \leq 1.5n+8t-6$.
\end{conjecture}
\begin{conjecture}
If $G$ is a graph with diameter $d$, then $\pi_t(G) \le \left(
    \frac{2^{\left \lceil \frac{d}{2} \right \rceil} - 1}{\left \lceil
    \frac{d}{2} \right \rceil} \right) n + 2^d (t-1) + f(d)$,
for some function $f$ that depends on $d$ only.
\label{t pebbling general diam}
\end{conjecture}
We show Conjecture~\ref{diamconj} holds for sufficiently large
$t$ after giving one lemma.
\begin{lemma}
\label{pi_t >= 2^d t}
For every graph $G$ and every integer $t \geq 1$, we have $\pi_t(G) \geq
2^{\diam(G)} t$.
\end{lemma}
\textbf{Proof}: We simply note that placing $2^{\diam(G)}t - 1$
pebbles on some vertex $v$ would create a situation from which we
could not move $t$ pebbles onto another vertex whose distance from $v$
is $\diam(G)$.\Bx
\begin{theorem}
\label{diambound}
For every graph $G$ with $n$ vertices, and for every $t \geq
\left\lceil \frac{n-1}{\diam(G)} \right\rceil$, we have $\pi_{t+1}(G)
\leq \pi_t(G) + 2^{\diam(G)}$.
\end{theorem}
\textbf{Proof}: We let $d = \diam(G)$.  By Lemma~\ref{pi_t >= 2^d t}
we have $\pi_t(G) \geq 2^d t$, so
\[
\pi_t(G) + 2^d \geq 2^d (t+1) \geq 2^d \left( \frac{n-1}{d} + 1
\right) = \frac{2^d}{d} (n-1) + 2^d \geq \frac{2^d - 1}{d} (n-1) + 1.
\]
Therefore, by Theorem~\ref{radius bound}, if we have $\pi_t(G) + 2^d$
pebbles on $G$, putting the first pebble on any target vertex costs at
most $2^d$ pebbles, so we can use the remaining $\pi_t(G)$ pebbles to
put $t$ additional pebbles on the target.~\Bx

In keeping consistent with the definitions of fractional analogues of other 
graph invariants, the fractional pebbling number was defined 
in~\cite{hurlbert URL} as follows. \\
\textbf{Definition} (Hurlbert~\cite{hurlbert URL}): The 
\emph{fractional pebbling number}
$\hat{\pi}(G)$ is given by $\displaystyle{ \hat{\pi}(G) = \liminf_{t
\rightarrow \infty} \frac{\pi_{t}(G)}{t}.}$ \\
In~\cite{ScUl}, we find a similar form for the definitions of the fractional 
analogues of chromatic number, clique number, matching number, and others.  
We use Theorem~\ref{diambound} to prove that $\hat{\pi}(G) =
2^{\diam(G)}$ for every graph $G$, as conjectured in~\cite{hurlbert
URL}.
\begin{theorem}
For every graph $G$, we have $\hat{\pi}(G) = 2^{\diam(G)}$.
\label{fractional pebbling}
\end{theorem}
\textbf{Proof}: We let $s = \left\lceil \frac{n-1}{\diam(G)}
\right\rceil$.  Applying Theorem~\ref{diambound} inductively on $t$
gives $\pi_t(G) \leq \pi_s(G) + (t-s) 2^{\diam(G)}$ for all $t \geq
s$.  Given $\epsilon > 0$, we let $x = \pi_s(G) - 2^{\diam(G)} s \geq
0$.  Then for any $t \geq \max(\frac{x}{\epsilon}, s)$ we have
\[
2^{\diam(G)} t \leq \pi_t(G) \leq \pi_s(G) + (t-s) 2^{\diam(G)} =
2^{\diam(G)} t + x.
\]
Dividing by $t$ gives
\[
2^{\diam(G)} \leq \frac{\pi_t(G)}{t} \leq 2^{\diam(G)} + \frac{x}{t} \leq
2^{\diam(G)} + \epsilon.
\]
Thus, $\hat{\pi}(G) = \displaystyle{\liminf_{t \rightarrow \infty}
\frac{\pi_{t}(G)}{t}} = 2^{\diam(G)}$.\Bx

\subsection{Optimal Fractional Pebbling Numbers}\label{OFP}

In this section, we see that
optimal pebbling can be modeled nicely as an optimization problem.
This in turn leads to a nice combinatorial interpretation of the
optimal fractional pebbling number of a graph.  We use this
interpretation to obtain a resulting property of vertex-transitive
graphs.  We begin by giving the generalization of a distribution to
allow non-integral amounts of pebbles to be placed on each vertex. \\
\textbf{Definition} (Moews~\cite{optimal hypercubes}): For a graph $G$, 
a function $D : V \rightarrow
\mathbb{R}^{\geq 0}$ is called a \emph{continuous distribution} on
$G$.  As in an integer-valued distribution, the \emph{size} of $D$ is
given by $|D|=\displaystyle{\sum_{v\in V}D(v)}$. \\
We give the following definition, which serves to generalize the
notion of a pebbling move. \\
\textbf{Definition} (Moews~\cite{optimal hypercubes}): 
A \emph{continuous pebbling move} of size
$\alpha\in\mathbb{R}^+$ from a vertex $v$, which has at least
$2\alpha$ pebbles, to a vertex $u\in N(v)$ removes $2\alpha$ pebbles
from $v$ and places $\alpha$ pebbles on $u$. \\
Thus, the pebbles are no longer discrete objects.  Instead, they can
be viewed as infinitely divisible ``piles''.  Nevertheless, for a
vertex $v$, a continuous distribution $D$, and a nonnegative real
number $\alpha$, if $D(v)=\alpha$, then we say that there are $\alpha$
pebbles on $v$ under $D$. \\
\textbf{Definition}: A continuous distribution $D$ on a graph $G$ is
called \emph{optimal} if the following two conditions hold.
\begin{enumerate}
\item For every $v\in V$, one pebble can be moved to $v$ after some
sequence of continuous pebbling moves, starting from $D$.
\item If $D'$ is a continuous distribution on $G$ with $|D'|<|D|$,
 then there is some $v\in V$ which cannot be reached with one pebble
 after any sequence of continuous pebbling moves, starting from $D'$.
\end{enumerate}

Recall that for a graph $G$ and an integer $t \geq 1$, $\pi_t^*(G)$ is
the size of the smallest $t$-fold solvable distribution of pebbles on
$G$.  Thus, given a $t$-fold solvable distribution $D$ on $G$, every
$v \in V$ must have a corresponding sequence of pebbling moves that
places $t$ pebbles on $v$, starting from $D$.  Let
$V=\{v_1,v_2,\ldots,v_n\}$.  For all $i$, $j$, and $k$, let
$p_{i,j,k}$ denote the number of pebbling moves from $v_j$ to $v_k$
in the sequence of moves which places a pebble on $v_i$.  Let us refer
to the following integer optimization problem as \textsf{OPT}. \\
\textbf{The \textsf{OPT} Integer Optimization Problem}: Minimize
$\displaystyle{\sum_{i=1}^n D(v_i)}$ subject to the following
constraints for each $i, j,$ and $k$ with $1 \leq i,j,k \leq n$:
\[ \begin{array}{c}
D(v_i) + \displaystyle{\sum_{v_j \sim v_i} (p_{i,j,i}-2p_{i,i,j})}
\geq t \\
D(v_k) + \displaystyle{\sum_{v_j \sim v_k} (p_{i,j,k}-2p_{i,k,j})} \geq 0 \\
D(v_i) \in \mathbb{N} \\
p_{i,j,k} \in \mathbb{N} \\
\end{array} \]
Clearly, every $t$-fold solvable distribution $D$ on $G$ results in a feasible 
solution to \textsf{OPT}.  Indeed, every pebbling move from a vertex removes 
two pebbles from it and every pebbling move to a vertex adds a pebble to it.  
Thus, after any sequence of pebbling moves which places at least $t$ pebbles 
on vertex $v_i$, every vertex must end up with a nonnegative number of 
pebbles while $v_i$ ends up with at least $t$ pebbles.   Conversely, Watson
\cite{Wats} shows that every feasible solution to \textsf{OPT} results in
a $t$-fold solvable distribution on $G$.  Thus, the solution to
\textsf{OPT} is equal to $\pi_t^*(G)$. \\
We give the following definition, which is similar to that of
$\hat{\pi}(G)$. \\
\textbf{Definition}: The \emph{optimal fractional pebbling number} 
$\hat{\pi}^*(G)$ is given by $\displaystyle{\hat{\pi}^*(G) = 
\liminf_{t \rightarrow \infty} \frac{\pi_{t}^*(G)}{t}.}$

Suppose that we desire a combinatorial interpretation for
$\hat{\pi}^*(G)$.  In this spirit, suppose we relax the integer
constraints in \textsf{OPT} and set $t=1$. Let us refer to the following
optimization problem as \textsf{FRAC OPT}. We denote its solution
$\ofpn(G)$, as in ~\cite{optimal hypercubes}, where this quantity is referred to
as the \emph{continuous optimal pebbling number of
$G$}. \\
\textbf{The \textsf{FRAC OPT} Optimization Problem}: Minimize
$\displaystyle{\sum_{i=1}^n D(v_i)}$ subject to the following
constraints for each $i, j,$ and $k$ with $1 \leq i,j,k \leq n$:
\[ \begin{array}{c}
D(v_i) + \displaystyle{\sum_{v_j \sim v_i} (p_{i,j,i}-2p_{i,i,j})}
\geq 1 \\
D(v_k) + \displaystyle{\sum_{v_j \sim v_k} (p_{i,j,k}-2p_{i,k,j})} \geq 0 \\
D(v_i) \geq 0 \\
p_{i,j,k} \geq 0 \\
\end{array} \]

We show that $\ofpn(G)$ is equal to the optimal fractional
pebbling number of $G$.

\begin{fact}\label{fracoptpebbthm}
For every graph $G$, $\ofpn(G) = \hat{\pi}^*(G)$.
\end{fact}
\textbf{Proof}: Let $G$ be a graph, with $V = \{ v_1, v_2, \ldots, v_n
\}$.  We first show $\ofpn(G) \leq \hat{\pi}^*(G)$.  For an
integer $t \geq 1$, let $D$ be a $t$-fold solvable distribution on
$G$ with $|D|=\pi_t^*(G)$.  Then, for every $v_i \in V$, there are
$D(v_i)$ pebbles initially on $v_i$ and there is some sequence of
pebbling moves which places $t$ pebbles on $v_i$.  This gives a
solution to \textsf{OPT}.  In this solution, let $p_{i,j,k}$ be
defined as above.  Now, let $D'(v_i) = \frac{D(v_i)}{t}$ and let
$p'_{i,j,k} = \frac{p_{i,j,k}}{t}$ for all $i$, $j$, and $k$.
This gives a feasible solution to \textsf{FRAC OPT} with $|D'| =
\frac{\pi_t^*(G)}{t}$.  This solution may or may not be optimal.
Since this holds for any integer $t \geq 1$, we have
$\ofpn(G) \leq \hat{\pi}^*(G)$.

We now show that $\ofpn(G) \geq \hat{\pi}^*(G)$.  Suppose
we have a feasible solution to \textsf{FRAC OPT}, with values denoted
$\overline{D}(v_i)$ and $\overline{p}_{i,j,k}$ for all $i$, $j$,
and $k$.  We may assume that every $\overline{D}(v_i)$ and
$\overline{p}_{i,j,k}$ is rational, since all of the coefficients
are integers.  Let $t$ be the least common multiple of the
denominators of these values.  Let $D'(v_i)=t\overline{D}(v_i)$ and
let $p'_{i,j,k}=t\overline{p}_{i,j,k}$ for all $i$, $j$, and
$k$.  This gives a feasible solution to \textsf{OPT} and thus a $t$-fold
solvable distribution $D'$ on $G$.  Clearly, $\frac{|D'|}{t}$ is the
value of the rational solution we were given.  However, $D'$ may not
be the smallest $t$-fold solvable distribution on $G$.  Furthermore,
we can let $D''(v_i)=ts\overline{D}(v_i)$ and let
$p''_{i,j,k}=ts\overline{p}_{i,j,k}$ for all $i$, $j$, and $k$
for any positive integer $s$ to obtain a $ts$-fold solvable
distribution on $G$.  Thus, $\ofpn(G) \geq
\hat{\pi}^*(G)$.\Bx

The following corollary provides a combinatorial interpretation for
$\hat{\pi}^*(G)$.
\begin{corollary}
The size of an optimal continuous distribution on a graph $G$ is equal
to $\hat{\pi}^*(G)$.
\end{corollary}
\textbf{Proof:} From the definition, we see that the size of an
optimal continuous distribution on a graph $G$ is equal to the
solution to the optimization problem \textsf{FRAC OPT}.  The result
follows from Fact \ref{fracoptpebbthm}.\Bx

In Theorem \ref{verttransthm}, we show that every vertex-transitive graph
has an optimal continuous distribution which is uniform.  We start with
some lemmas, beginning with the following self-evident weight argument.
\begin{lemma}
\label{weightarg}
Let $D$ be a continuous distribution on a graph $G$.  Then there is a
sequence of continuous pebbling moves starting from $D$ which places a
pebble on $\br\in V$ if and only if $\displaystyle{\sum_{v\in V} D(v)
2^{-\dist(v,\br)} \geq 1}$.~\Bx
\end{lemma}
The following lemma is obvious, but useful.
\begin{lemma}
\label{constant function}
If $G=(V,E)$ is a vertex-transitive graph, then the function $f : V
\rightarrow \mathbb{R}^+$ given by $f(u) = \displaystyle{\sum_{v \in
V} 2^{-\dist(v, u)}}$ is constant for all $u$.~\Bx
\end{lemma}
\begin{lemma}
\label{fracoptlem}
If $D$ and $D'$ are continuous distributions on a vertex-transitive
graph $G$ and 
\begin{equation}
\label{losing weight}
\sum_{u \in V} D(u) 2^{-\dist(v,u)} \leq \sum_{u \in V} D'(u)
2^{-\dist(v,u)}
\end{equation}
for all $v \in V$, then $|D| \leq |D'|$.
\end{lemma}
\textbf{Proof}: Let $G=(V,E)$ be a vertex-transitive graph.  Summing
both sides of~(\ref{losing weight}) over all $v \in V$, we find
\[
\sum_{v \in V} \sum_{u \in V} D(u) 2^{-\dist(v,u)} \leq \sum_{v \in
V}\sum_{u \in V} D'(u) 2^{-\dist(v,u)}
\]
Switching the order of the summation gives us
\[
\sum_{u \in V} D(u) \sum_{v \in V} 2^{-\dist(v,u)} \leq \sum_{u \in
V} D'(u) \sum_{v \in V} 2^{-\dist(v,u)}.
\]
But by Lemma~\ref{constant function}, $\displaystyle{\sum_{v \in V}
2^{-\dist(v,u)}}$ is a constant for all $u \in V$, so dividing by
this constant gives us $\displaystyle{\sum_{u \in V} D(u) \leq \sum_{u
\in V} D'(u)}$, or $|D| < |D'|$. \Bx

We are now ready to show the main result for this section.

\begin{theorem}
\label{verttransthm}
If $G$ is a vertex-transitive graph, an optimal continuous distribution
on $G$ is obtained by putting $\frac{1}{m}$ pebbles on each vertex in
$G$, where $m$ is the constant $\displaystyle{\sum_{v\in V}
2^{-\dist(v,u)}}$ from Lemma~\ref{constant function}.  Therefore,
$\hat{\pi}^*(G) = \frac{n}{m}$.
\end{theorem}
\textbf{Proof}: Let $D$ be the distribution in question.  Note that
for any root $\br\in V$, we have
\[
\sum_{v \in V} D(v) 2^{-\dist(v,\br)} = \frac{1}{m} \sum_{v\in V}
2^{-\dist(v,\br)} = 1,
\]
so by Lemma~\ref{weightarg}, starting from $D$, the root $\br$ can
receive a pebble by making continuous pebbling moves toward $\br$.
Therefore, $\hat{\pi}^*(G)\leq|D|$.

Now let $D'$ be another continuous distribution from which one pebble
can be moved to $\br$.  By Lemma~\ref{weightarg}, we have
\[
\sum_{v\in V} D'(v) 2^{-\dist(v,\br)} \geq 1 = \sum_{v\in V} D(v)
2^{-\dist(v,\br)}
\]
for all $v \in V$, and by Lemma~\ref{fracoptlem}, this implies $|D'|
\geq |D|$.  Therefore, $D$ is optimal, so $\hat{\pi}^*(G)=|D|$.\Bx
\\
Corollary~\ref{vertex transitive graphs} gives $\hat{\pi}^* (G)$ for
several vertex-transitive graphs.  Moews~\cite{optimal hypercubes}
also proved part~\ref{hypercube}.
\begin{corollary}
Let $k$ and $n$ be positive integers.  Then we have the following.
\begin{enumerate}
\item $\displaystyle{\hat{\pi}^*(Q^k)=\left(\frac{4}{3}\right)^k}$
where $Q^k$ denotes the $k$-dimensional hypercube.
\label{hypercube}
\item $\displaystyle{\hat{\pi}^*(K_n) = \frac{2n}{n+1}}$.
\item If $k \geq 2$, then $\displaystyle{\hat{\pi}^*(C_{2k}) =
  \frac{k2^{k+1}}{3(2^k-1)}}$.
\item $\displaystyle{\hat{\pi}^* (C_{2k+1}) =
 \frac{(2k+1)(2^{k-1})}{3(2^{k-1})-1}}$.
\end{enumerate}
\label{vertex transitive graphs}
\end{corollary}
\textbf{Proof}: By Theorem~\ref{verttransthm}, in each case it
suffices to find the value of $m$.  For the hypercube, if we fix a
target $\br$, there are ${k \choose i}$ vertices whose distance from $\br$
is $i$.  We compute $m$ as follows, using the Binomial Theorem.
\[
m = \sum_{v\in V} 2^{-\dist(v,\br)} = \sum_{i=1}^k {k \choose i}
\frac{1}{2^i} = \left( \frac{3}{2} \right)^k.
\]
Therefore, $\hat{\pi}^* (Q^k) = \frac{n}{m} = \frac{2^k}{\left(
  \frac{3}{2} \right)^k} = \left( \frac{4}{3} \right)^k$.

For $K_n$ every vertex $v \neq r$ has $\dist(v, \br) = 1$, so 
$$m = 1 +
\sum_{v \in V; v \neq \br} \frac{1}{2} = 1 +
\frac{n-1}{2} = \frac{n+1}{2},$$ and $\hat{\pi}^* (K_n) =
\frac{n}{\frac{n+1}{2}} = \frac{2n}{n+1}$.

For $C_n$ we assume the vertex set is $\{ x_0, x_1, \ldots, x_{n-1}
\}$ and that $\br = x_0$ is the target.  If $n = 2k$, we let $A = \{ x_i
: i < k \}$, and we note that for every $x_{k+i}$ with $0 \leq i \leq
k-1$ we have $\dist(x_{k+i}, x_0) = k-i$.  Therefore, computing $m$
gives
\[
m = \sum_{v\in V} 2^{-\dist(v,x_0)} = \sum_{v\in A} 2^{-\dist(v,x_0)}
+ \sum_{v \not \in A} 2^{-\dist(v,x_0)} = \sum_{i=0}^{k-1} 2^{-i} +
\sum_{i=0}^{k-1} 2^{-(k-i)}.
\]
Substituting $j = k-1-i$ in the last summation gives
\[
m = \sum_{i=0}^{k-1} 2^{-i} + \sum_{j=0}^{k-1} 2^{-(j+1)} =
\sum_{i=0}^{k-1} 2^{-i} + \frac{1}{2} \sum_{j=0}^{k-1} 2^{-j} =
\frac{3}{2} \left( 2 - \frac{1}{2^{k-1}} \right) = \frac{3(2^k -
1)}{2^k}.
\]
Therefore, $\hat{\pi}^* (C_{2k}) = \frac{n}{m} = \frac{2k (2^k)}{3(2^k
  - 1)} = \frac{k2^{k+1}}{3(2^k-1)}$. 

Finally, for $C_{2k+1}$ we let $A = \{ x_i : 1 \leq i \leq k \}$ and
$B = \{ x_i : k+1 \leq i \leq 2k \}$.  Now $\dist(x_{k+i}, x_0) =
k-i+1$, so we have
\[
m = \sum_{v \in V} 2^{-\dist(v, x_0)} = 2^{-\dist(x_0, x_0)} + \sum_{v
\in A} 2^{-\dist(v, x_0)} + \sum_{v \in B} 2^{-\dist(v, x_0)} = 1 +
\sum_{i=1}^k 2^{-i} + \sum_{i=1}^{k} 2^{-(k-i+1)}.
\]
Now substituting $j = k-i+1$ gives
\[
m = 1 + \sum_{i=1}^k 2^{-i} + \sum_{j=1}^{k} 2^{-j} = 1 + 2
\sum_{i=1}^k 2^{-i} = 1 + 2 \left( 1 - \frac{1}{2^k} \right) = 3 -
\frac{1}{2^{k-1}} = \frac{3 (2^{k-1}) - 1}{2^{k-1}}.
\]
Therefore, $\hat{\pi}^* (C_{2k+1}) = \frac{n}{m} =
\frac{(2k+1)(2^{k-1})}{3(2^{k-1})-1}$.\Bx

\bibliographystyle{plain}

\end{document}